\documentclass[11pt]{amsart}
\usepackage{a4}
\usepackage{amsfonts}
\usepackage{amsmath}
\usepackage{amssymb}

\usepackage{hyperref}
\usepackage{amscd}
\date{}
\numberwithin{equation}{section}

\swapnumbers
\theoremstyle{plain}
\newtheorem{Lemma}{Lemma}
\numberwithin{Lemma}{section}

\newtheorem{Corollary}[Lemma]{Corollary}

\newtheorem{Theorem}[Lemma]{Theorem}

\theoremstyle{definition}

\theoremstyle{remark}
\newtheorem{Remark}[Lemma]{Remark}

\newtheorem{Example}[Lemma]{Example}

\newcommand\A{\mathbb A}
\newcommand\B{\mathbb B}

\newcommand\C{\mathbb C}
\newcommand\HH{\mathbb H}

\newcommand\Q{\mathbb Q}
\newcommand\R{\mathbb R}
\newcommand\Z{\mathbb Z}

\renewcommand\Re{\operatorname{Re}}
\renewcommand\Im{\operatorname{Im}}
\newcommand\sign{\operatorname{sign}}

\newcommand\id{\mathrm{id}}
\newcommand\im{\operatorname{im}}
\newcommand\ind{\operatorname{ind}}
\newcommand\Ind{\operatorname{Ind}}
\newcommand\punkt{\mathord{\,\cdot\,}}
\newcommand\ch{\mathrm{ch}}
\newcommand\tr{\operatorname{tr}}

\newcommand\rk{\operatorname{rk}}
\newcommand\ev{{\mathrm{ev}}}
\newcommand\odd{{\mathrm{odd}}}
\newcommand\SF{\operatorname{sf}}
\newcommand\Spec{\operatorname{Spec}}
\newcommand\lk{\operatorname{lk}}

\newcommand\eps{\varepsilon}

\begin{document}

\title[Eta invariants]{Computations and Applications of $\eta$ invariants}
\author{Sebastian Goette}
\address{Mathematisches Institut\\
Universit\"at Freiburg\\
Eckerstr.~1\\
79104 Freiburg\\
Germany}
\email{sebastian.goette@math.uni-freiburg.de}
\thanks{Supported in part by DFG special programme
``Global Differential Geometry''}
\begin{abstract}
  We give a survey on $\eta$-invariants including methods of computation
  and applications in differential topology.
\end{abstract}

\maketitle

\section*{Introduction}
The $\eta$-invariant has been introduced by Atiyah, Patodi and Singer as a boundary contribution in an index theorem for elliptic operators in the series of papers~\cite{APS0}--\cite{APS3}. There are several invariants of odd-dimensional manifolds~$M$ in differential topology that are originally defined by finding a compact manifold~$N$ with boundary~$\partial N=M$ and evaluating certain characteristic numbers on~$N$. The Atiyah-Patodi-Singer index theorem~\ref{A1.T1} often allows to compute these invariants in terms of $\eta$-invariants and other magnitudes that can be defined directly on~$M$ without choosing~$N$ first. Sometimes this leads to generalisations of these invariants to manifolds that are not 0-cobordant. However, to determine such an invariant for a given manifold~$M$, one needs ways to compute $\eta$-invariants of operators defined on~$M$ without using the Atiyah-Patodi-Singer index theorem.

In the present paper, we give a short survey on applications of $\eta$-invariants, with a focus on situations where the corresponding $\eta$-invariants can be computed. The $\eta$-invariant also appears in possible generalisations of the analytic torsion, in conformal geometry, and in the definition of certain smooth extensions of topological $K$-theory. To keep this article reasonably short, we will not touch upon these and several other issues.

We start by reviewing the definition of $\eta$-invariants as spectral invariants in section~\ref{A}. We also review the Atiyah-Patodi-Singer theorem and some of its immediate consequences. We list some examples where $\eta$-invariants have been computed directly. Sometimes it is easier to compute $\eta$-invariants of modified operators first and to determine their difference to the original $\eta$-invariants, see sections~\ref{B3} and~\ref{C2}.

The Atiyah-Patodi-Singer theorem has generalisations in different directions. In section~\ref{B}, we consider families of manifolds, group actions and orbifolds. The corresponding generalisations of Theorem~\ref{A1.T1} involve generalisations of $\eta$-invariants that are sometimes easier to compute. 

In section~\ref{P}, we discuss the behaviour of $\eta$-invariants for direct images under proper maps and under gluing constructions. These methods sometimes give rise to explicit computations, see section~\ref{C5}.

Finally, in section~\ref{C}, we discuss some applications of $\eta$-invariants mainly to differential topology, but also to questions ranging from algebraic $K$-theory to Riemannian manifolds of positive scalar curvature.

The author wishes to thank M. Braverman, U. Bunke, D. Crowley, P. Piazza and A. Ranicki for some helpful comments and explanations.

\section{The Atiyah-Patodi-Singer \texorpdfstring{$\eta$}{eta}-invariant and related invariants}\label{A}

The~$\eta$-invariant of a selfadjoint elliptic differential operator on an odd-dimensional manifold~$M$ first appeared in the Atiyah-Patodi-Singer index theorem for manifolds~$N$ with boundary~$M$, which was announced in~\cite{APS0} and proved in~\cite{APS1}. They already noted that the~$\eta$-invariant was related to other topological invariants known at the time.

\subsection{The index theorem for manifolds with boundary}\label{A1}

Let~$N$ be a compact Riemannian manifold with boundary~$M = \partial N$, and let~$A \colon \Gamma(E^+) \to \Gamma(E^-)$ be an elliptic differential operator on~$N$. Assume that a neighbourhood~$U$ of~$M$ in~$N$ is isometric to a product~$M \times [0,\eps)$, that~$\nu \colon E^+|_U \to E^-|_U$ is a vector bundle isomorphism, and that~$E^\pm|_U$ are identified with~$E^\pm|_M \times [0,\eps)$ in such a way that on~$U$,
\begin{equation}\label{A1.1}
A|_{\Gamma(E^+|_U)} = \nu \circ \Bigl(\frac{\partial}{\partial t} + B\Bigr)\;,
\end{equation}
where~$\frac{\partial}{\partial t}$ denotes differentiation in the direction of~$[0,\eps)$ and~$B$ is a selfadjoint elliptic differential operator acting on smooth sections of~$E^+|_{\partial N}\to M$.

A typical example consists of a Dirac operator~$A = D_N^+$ on an even-dimensional manifold. With
$$ \nu = c_N\Bigl(\frac{\partial}{\partial t}\Bigr) \colon E^+|_U \to E^-|_U\;,$$
one can construct a Clifford multiplication~$c_M$ of~$TM$ on~$E^+|_M$, such that
$$ c_N(v) = \nu \circ c_M(v) \qquad \text{for all}\qquad v \in TM\;.$$
Then~\eqref{A1.1} holds with~$B = D_M$ a Dirac operator on the odd-dimensional boundary~$M$.

There are topological obstructions against local elliptic boundary conditions for the operator~$A$ on~$N$ of~\eqref{A1.1}. However, there are elliptic spectral boundary conditions on each connected component of~$M$, inspired by replacing~$U$ by an infinite cylinder~$M \times (-\infty,\eps)$ and imposing~$L^2$-boundary conditions in the special case that~$B$ is invertible. Concretely, one restricts~$A$ and its adjoint~$A^*$ to
\begin{equation}\label{A1.2}
  \begin{aligned}
    \Gamma_<(E^+)
    &=\bigl\{\,\sigma \in \Gamma(E^+) \bigm| P_\ge(\sigma|_M) = 0\,\bigr\}\;,\\
    \Gamma_\le(E^-)
    &=\bigl\{\,\tau \in \Gamma(E^-) \bigm| P_<(\nu^{-1}\tau|_M) = 0\,\bigr\}\;,
  \end{aligned}
\end{equation}
where~$P_\ge$, $P_< \colon \Gamma(E^+|_M) \to \Gamma(E^+|_M)$ denote spectral projections onto the nonnegative and the negative eigenspaces of~$B$, respectively. These are the so-called {\em APS boundary conditions}, and one defines
$$\ind_{\mathrm{APS}}(A)=\ker\bigl(A|_{\Gamma_<(E^+)}\bigr)-\ker\bigl(A^*|_{\Gamma_\le(E^-)}\bigr)\;.$$

The index of a suitable double of~$A$ is given as the integral of a local index density~$\alpha_0$ over the double of~$N$. In the case of a Dirac operator, we know from the Atiyah-Singer index theorem that
$$ \alpha_0 = \bigl(\hat{A}(TN,\nabla^{TN}) \wedge \ch(E/S,\nabla^{E})\bigr)^{\max}\;,$$
where~$\ch(E/S,\nabla^E)$ denotes the twist Chern character form,
see~\cite{BGV}.

The~$\eta$-invariant~$\eta(B)$ of~$B$ is defined as the value at~$s = 0$ of the meromorphic continuation of the~$\eta$-function, that is for~$\Re s\gg 0$ given by
\begin{equation}\label{A1.3}
\eta_B(s) = \sum_{\lambda \in \Spec(B)\backslash \{0\}}\sign(\lambda)\cdot|\lambda|^{-s} = \frac{1}{\Gamma\bigl(\frac{s+1}{2}\bigr)} \int_0^\infty t^{\frac{s-1}{2}}\tr(Be^{-tB^2})\, dt\;.
\end{equation}
It is proved in~\cite{APS3} that the~$\eta$-function indeed has a meromorphic continuation and that~$\eta(B)=\eta_B(0)$ is finite. For a Dirac operator~$B = D_M$, one can show directly that the integral expression in~\eqref{A1.3} converges for~$\Re s > -1$, see~\cite{BF}. One also defines
$$ h(B) = \dim \ker B\;.$$

\begin{Theorem}[Atiyah-Patodi-Singer, \cite{APS0}, \cite{APS1}]\label{A1.T1}
Let~$A$ be an elliptic differential operator on a compact manifold~$N$ with boundary~$M = \partial N$ as in~\eqref{A1.1}. Then the Fredholm index of~$A$ under the APS boundary conditions~\eqref{A1.2} is given as
$$ \ind_{\mathrm{APS}}(A) = \int_N \alpha_0 - \frac{\eta + h}{2}(B)\;.$$
\end{Theorem}

The signature operator of a~$4k$-dimensional compact oriented manifold with boundary is an important special case. Here, one considers the symmetric bilinear form on
$$ \im\bigl(H^{2k}(N,\partial N;\R) \to H^{2k}(N;\R)\bigr)$$
given by the evaluation of the cup product of two such classes on the relative fundamental cycle~$[N,\partial N]$. The signature of this form is denoted by~$\sign(N)$. 

On the odd-dimensional manifold~$M = \partial N$, the bundle~$\Lambda^\ev T^*M$ of even differential forms constitutes a Dirac bundle. The Hodge star operator~$*$ interchanges even and odd forms. The Dirac operator~$B = \pm (*d-d*)$ on~$\Lambda^\ev T^*M$ is usually named the {\em odd signature operator} on~$M$. 

\begin{Theorem}[Atiyah-Patodi-Singer, \cite{APS1}]\label{A1.T2} Let~$N$ be a~$4k$-dimensional manifold with totally geodesic boundary~$M = \partial N$, and let~$B$ denote the odd signature operator on~$M$, then
$$ \sign(N) = \int_N L(TN) - \eta(B)\;.$$
\end{Theorem}

Comparing with Theorem~\eqref{A1.T1}, one notes that the boundary operator consists of two copies of the odd signature operator. Also, the index of the signature operator on~$N$ under APS boundary conditions is not precisely~$\sign(N)$ due to the asymmetric treatment of~$\ker(B)$. In fact, if~$h(B)$ was present on the right hand side in Theorem~\eqref{A1.T2}, the equation would not be compatible with a change of orientation. The most prominent feature for applications is the fact that the signature of~$N$ is a topological invariant, in contrast to most other APS indices, which depend on the geometry of~$N$ near its boundary.

Some elementary properties of~$\eta$-invariants can be deduced directly from Theorem~\ref{A1.T1} and~\ref{A1.T2}. For simplicity, we will stick to Dirac operators, and we let~$B$ denote the odd signature operator.

If~$P(V,\nabla^V) \in \Omega^\bullet(M)$ denotes a Chern-Weil form associated to a vector bundle~$V \to M$ with connection~$\nabla^V$ and an invariant polynomial~$P$, we let~$\tilde{P}(V,\nabla^{V,0}\nabla^{V,1}) \in \Omega^\bullet(M)/\im d$ denote the Chern-Simons class satisfying
\begin{equation}\label{A1.4}
d\tilde{P}\bigl(V,\nabla^{V,0},\nabla^{V,1}\bigr) = P\bigl(V,\nabla^{V,1}\bigr) - P\bigl(V,\nabla^{V,0}\bigr)\;.
\end{equation}
The Dirac operator~$D$ on a Dirac bundle~$E \to M$ depends on smoothly on the Riemannian metric~$g$ on~$M$ and on a Clifford multiplication and a connection~$\nabla^E$ on~$E$ that is compatible with the Levi-Civita connection. Applying Theorems~\ref{A1.T1} and~\ref{A1.T2} to an adapted Dirac operator~$D_N$ on the cylinder~$N = M \times [0,1]$ gives a variation formula.

\begin{Corollary}[Atiyah-Patodi-Singer, \cite{APS1}]\label{A1.C3}
Let~$(g_s)_{s \in [0,1]}$ be a family of Rie\-mann\-ian metrics on~$M$ with Levi-Civita connections~$\nabla^{TM,s}$, and let~$(E,c_s,\penalty0\nabla^{E,s})_{s \in [0,1]}$ be compatible bundles with Dirac operators~$D_M^s$. Then
\begin{multline*}
\frac{\eta + h}{2}(D_M^1) - \frac{\eta + h}{2}(D_M^0) \\
= \int_M \Bigl(\tilde{\hat A}\bigl(TM,\nabla^{TM,0},\nabla^{TM,1}\bigr)\, \ch\bigl(E/S,\nabla^{E,0}\bigr) \\
-\hat{A}\bigl(TM,\nabla^{TM,1}\bigr)\, \tilde\ch\bigl(E/S,\nabla^{E,0},\nabla^{E,1}\bigr)\Bigr)\in\R/\Z\;. \tag{1}\label{A1.C3.1}
\end{multline*}
For the odd signature operator~$(B^s)_{s \in [0,1]}$, one has
\begin{gather*}
\eta(B^1) - \eta(B^0) = \int_M \tilde{L}\bigl(TM,\nabla^{TM,0},\nabla^{TM,1}\bigr)\in\R\;.
\tag{2}\label{A1.C3.2}
\end{gather*}
\end{Corollary}

Thus, $\eta$-invariants have similar variation formulas as Cheeger-Simons numbers, by which we mean the evaluation of (products of) Cheeger-Simons classes on the fundamental cycle of an odd-dimensional compact oriented manifold. We may think of Cheeger-Simons numbers as geometric~$\R/\Z$-valued refinements of integral characteristic classes of vector bundles, whereas~$\eta$-invariants are geometric~$\R/\Z$-valued refinements of indices of Dirac operators. In general, these numbers are difficult to compare.

\begin{Example}\label{A1.E4} Let~$M$ be a oriented three-manifold. Then the variation formula for~$\eta(B)$ of Corollary~\ref{A1.C3} becomes
$$ \eta(B^1) - \eta(B^0) = \frac{1}{3} \int_M \tilde{p}_1 \bigl(TM,\nabla^{TM,0},\nabla^{TM,1}\bigr)\;,$$
where~$\tilde{p}_1$ is the Chern-Simons class associated to the first Pontrijagin class. The variation formula for the associated Cheeger-Simons character~$\hat{p}_1$ of a general vector bundle is
$$ \Bigl(\hat{p}_1\bigl(E,\nabla^{E,1}\bigr) - \hat{p}_1\bigl(E,\nabla^{E,0}\bigr)\Bigr)[M] = \int_M \tilde{p}_1\bigl(E,\nabla^{E,0},\nabla^{E,1}\bigr)\;.$$
This implies that~$3\eta(B)$ is an~$\R$-valued refinement of~$\hat{p}_1(TM,\nabla^{TM})[M]$. For other vector bundles, we do not get a natural~$\R$-valued refinement of~$\hat{p}_1(E,\nabla^E)[M] \in \R/\Z$ due to the presence of a nontrivial gauge group.
\end{Example}

For higher dimensional manifolds, the situation is more complicated due to the formulas for the multiplicative sequences,
\begin{align}
\hat{A} &= 1 - \frac{p_1}{24} + \frac{7p_1^2 - 4p_2}{2^7 \cdot 3^2 \cdot 5} - \frac{31p_1^3 - 44p_1p_2 + 16p_3}{2^{10} \cdot 3^3 \cdot 5 \cdot 7} \pm \dots\;, \label{A1.5} \\
L &= 1 + \frac{p_1}{3} - \frac{p_1^2 - 7p_2}{45} + \frac{2p_1^3 - 13p_1p_2 + 62p_3}{3^3 \cdot 5 \cdot 7} \pm \dots \label{A1.6}
\end{align}

Finally, one can use Corollary~\ref{A1.C3} to count sign changes of eigenvalues of~$D_M^s$ for~$s \in [0,1]$ by comparing the actual difference of~$\eta$-invariants with the value predicted by the local Chern-Simons variation terms. The so-called {\em spectral flow} of the family~$(D_M^s)_{s \in [0,1]}$ is given by
\begin{multline}\label{A1.7}
\SF\bigl((D^s_M)_{s \in [0,1]}\bigr) = \frac{\eta + h}{2}(D_M^1) - \frac{\eta + h}{2}(D_M^0) \\
-\int_M\Bigl(\tilde{\hat A}\bigl(TM,\nabla^{TM,0},\nabla^{TM,1}\bigr)\ch\bigl(E/S,\nabla^{E,0}\bigr)\\ + \hat{A}\bigl(TM,\nabla^{TM,1}\bigr)\tilde{\ch}\bigl(E/S,\nabla^{E,0},\nabla^{E,1}\bigr)\Bigr) \in \Z\;.
\end{multline}
If the Dirac bundles for~$s = 0$ and~$s = 1$ are isomorphic, then the spectral flow defines an odd index~$\SF\colon K^1(M)\to\Z$, see~\cite{APS3}.

\subsection{Direct computation of \texorpdfstring{$\eta$}{eta}-invariants}\label{A2}
For generic Riemannian manifolds, it seems impossible to determine the spectrum of a given differential operator~$B$. And even if one succeeds, one often needs techniques from analytic number theory in order to describe the function~$\eta_B(s)$ explicitly and compute its special value~$\eta(B) = \eta_B(0)$ at~$s = 0$. This section is devoted to a few examples where this has been done. All examples are locally homogeneous spaces, and representation theory plays a prominent role in the determination of the relevant spectrum.

For the operator~$B_\lambda = i\frac{d}{dt} + \lambda$ on a circle of length~$2\pi$, the~$\eta$-invariant is computed in~\cite{APS0} as
$$ \eta(B_\lambda) = \begin{cases} 0 &\lambda \in \Z\;, \\ 1-2(\lambda - n) \qquad & \lambda \in (n,n+1)\;. \end{cases}$$

Next, consider three-dimensional Berger spheres. Thus, one rescales the fibres of the Hopf fibration~$S^3 \to S^2$ by~$\lambda > 0$, while the metric orthogonal to the fibres is unchanged. This metric is still~$U(2)$-invariant. Let~$D_\lambda$ denote the untwisted Dirac operator on~$S_\lambda^3$. Using~$U(2)$-invariance and a suitable Hilbert basis of sections of the spinor bundle, Hitchin determined the eigenvalues of~$D_\lambda$ in~\cite{Hi1} as follows,
\begin{align*}
& \frac{\lambda}{2} + \frac{p}{\lambda} & \qquad & \text{with multiplicity}\, 2p\;,\\
& \frac{\lambda}{2} \pm \frac{\sqrt{4pq \lambda^2 + (p-q)^2}}{\lambda} && \text{with multiplicity}\, p+q\;,
\end{align*}
for~$p$, $q > 0$. From these values, Hitchin computes the~$\eta$-invariant explicitly and obtains
$$ \eta(D_\lambda) = -\frac{(\lambda^2-1)^2}{6}$$
for~$0 < \lambda < 4$. For larger values of~$\lambda$, the formula holds only up to spectral flow, see~\eqref{A1.7}. In his diploma thesis~\cite{Habel}, Habel does analogous computations for a few higher dimensional Berger spheres. Bechtluft-Sachs gets the same result in low dimensions by applying Theorem~\ref{A1.T1} to disk bundles over~$\C P^n$ \cite{BS}.

Let~$\Gamma \subset PSL(2,\R)$ be a cocompact subgroup, then~$M = PSL(2,\R)/\Gamma$ is a compact three-manifold, Seifert fibred over a hyperbolic surface. Seade and Steer compute~$\eta(D_\lambda)$ when~$\Gamma$ is a Fuchsian group~\cite{SeSt}. The parameter~$\lambda$ refers to the length of the generic fibre as in the case of the Berger sphere. As in Hitchin's computations, representation theory plays a prominent role in these computations. The results of Seade and Steer are generalised to noncompact quotients of finite volume by Loya, Moroianu and Park in~\cite{LMP}.

The spectrum of an untwisted Dirac~$D$ operator on a flat torus~$M$ depends on the spin structure and is always symmetric, so~$\eta(D) = 0$. However, among the Bieberbach manifolds~$M/\Gamma$, where~$\Gamma \subset SO(n)$ is a finite subgroup that acts freely on~$M$, there exist examples with asymmetric spectra and non-vanishing~$\eta$-invariants. Pf\"affle computes the spectra and the~$\eta$-invariants of all three-dimensional examples in~\cite{Pf}. Higher-dimensional examples are studied by Sadowski and Szczep\'anski in~\cite{SS}, by Miatello and Podest\'a in~\cite{MiPo}, and by Gilkey, Miatello and Podest\'a in~\cite{GMP}. In all these cases, the $\eta$-invariant of an untwisted Dirac operator can be expressed in number theoretic terms.

Similar computations are also possible for spherical space forms~$S^n/\Gamma$ with~$\Gamma \subset SO(n)$ a finite subgroup.  In~\cite{Cm}, Cisneros-Molina gives a general formula for the spectra of Dirac operators on~$M = S^3/\Gamma$ twisted by flat vector bundles and computes the corresponding~$\eta$-invariants. These~$\eta$-invariants are closely related to the~$\xi$-invariant of~$M$ and the~$\Gamma$-equivariant~$\eta$-invariants of~$S^3$, see section~\ref{C1}. Seade~\cite{Seade} and Tsuboi~\cite{Tsuboi} compute~$\eta$-invariants for certain spherical space forms as average over equivariant~$\eta$-invariants, see also~\cite{Bspin}. Degeratu extends these computations to orbifold quotients in~\cite{Degeratu} and exhibits a relation with the Molien series.

In~\cite{Mill}, Millson expresses the $\eta$-invariant of the odd signature operator on a compact hyperbolic manifold as a special value of a $\zeta$-function associated to the closed geodesics on~$M$ and their holonomy representations and Poincar\'e maps. This result is extended to Dirac operators on locally symmetric spaces~$M$ of noncompact type by Moscovici-Stanton~\cite{MS}. A generalisation to the finite-volume case is given by Park~\cite{Park}.

The $\eta$-invariant of a Dirac operator on an interval~$[0,1]$,
twisted by a symplectic vector space~$(V,\omega)$,
with different Lagrangian boundary conditions~$L_0$, $L_1\subset\ker(D_N)$
is computed by Cappell, Lee and Miller in~\cite{CLM}.
It is related to the Maslov index of Lagrangians in~$(V,\omega)$.
Indeed, Maslov indices naturally occur when considering $\eta$-invariants
on manifolds with boundary,
for example in generalisations of Theorem~\ref{P2.T1}.

\section{Families, group actions, and orbifolds}\label{B}
Instead of regarding Dirac operators on a single manifold~$N$, one may consider families of manifolds, or manifolds with the action of some Lie-group, or even orbifolds with boundary. Under certain conditions, Theorems~\eqref{A1.T1} and~\eqref{A1.T2} extend to these situations. We state a few of these generalisations below and indicate relations between them.  We also explain how ordinary $\eta$-invariants can be computed using equivariant methods.

\subsection{Families of manifolds with boundary}\label{B1} 
Assume that~$p \colon W \to B$ is a proper submersion with typical even-dimensional fibre~$N$, such that the fibrewise boundaries form another submersion~$V \to B$. Let~$g^{TN}$ be a fibrewise Riemannian metric and let~$T^HW \to W$ be a horizontal complement for the fibrewise tangent bundle~$TN = \ker dp \subset TW$. These data uniquely define a generalised Levi-Civita connection~$\nabla^{TN}$ on~$TN \to W$. Let~$E = E^+ \oplus E^- \to W$ be a fibrewise Dirac bundle, i.e., $TN$ acts on~$E$ by Clifford multiplication, and there is a compatible metric~$g^E$ and a compatible connection~$\nabla^E$ on~$E$. Then one can define a family of Dirac operators~$D_N$ on the fibres of~$p$. We assume that condition~\eqref{A1.1} is satisfied on each fibre, and we also assume that the kernels of the boundary operators~$D_M$ form a family over~$B$. Then let us assume for simplicity that the kernels of the family~$D_N$ under APS boundary conditions also form a family over~$B$.

In this situation, there exist natural families of Bismut-Levi-Civita superconnections~$(\A_t)_{t \in (0,\infty)}$ and~$(\B_t)_{t \in (0,\infty)}$ on the infinite dimensional vector bundles~$p_*E \to B$ and~$p_*(E^+|_V) \to B$. These superconnections define ordinary connections~$\nabla^H$, $\nabla^{K^\pm}$ on~$H = \ker(D_M) \to B$ and on~$K^\pm = \ker(D_N^\pm) \to B$. The~$\eta$-invariant generalises to a natural~$\eta$-{\em form}
\begin{equation}\label{B1.1}
\tilde{\eta}(\B) = \frac{1}{\sqrt{\pi}}\int_0^\infty \tr\biggl(\frac{\partial \B_t}{\partial t}\, e^{-\B_t^2}\biggr) \, dt \in \Omega^\bullet(B)\;.
\end{equation}
Note that the component of degree 0 is~$\tilde{\eta}(\B)^{[0]} = \frac{\eta}{2}(B)$.

\begin{Theorem}[Bismut-Cheeger, \cite{BC1}, \cite{BC2}, \cite{BC3}]\label{B1.T1}
Under the assumptions above,
\begin{multline*}
\ch\bigl(K^+,\nabla^{K^+}\bigr) - \ch\bigl(K^-,\nabla^{K^-}\bigr) = \int_{W/B} \hat{A}\bigl(TN,\nabla^{TN}\bigr)\,\ch\bigl(E/S,\nabla^E\bigr) \\
-\tilde{\eta}(\B) - \frac{1}{2} \ch\bigl(H,\nabla^H\bigr) \in H^\bullet(B;\R)\;.
\end{multline*}
\end{Theorem}
There exists a similar generalisation of Theorem~\ref{A1.T2}. Note that the kernels of the signature operator and the odd signature operator on the boundaries automatically form bundles over~$B$ by Hodge theory.

Melrose and Piazza relax the condition that the kernels of the boundary operator~$D_M$ form a bundle over~$B$. For the definition of boundary conditions, it is sufficient to have a spectral section, see \cite{MP}. It is also not necessary to demand that~$\ker D_N$ forms a bundle over~$B$, since the virtual index bundle always exists.

\subsection{Group actions on manifolds with boundary}\label{B2}
Theorems~\ref{A1.T1} and~\ref{A1.T2} generalise to manifolds with group actions in the same way that the Atiyah-Singer index theorem becomes the Atiyah-Segal fixpoint theorem. In particular, from invariants on the boundary one can conclude the existence of fixpoints in the interior.

Let~$N$ and~$D_N$ be as in section~\ref{A1}. Let~$G$ be a group that acts on~$N$ by isometries. Assume that this action also lifts to the Dirac bundle~$E$, and that the induced action on~$\Gamma(E)$ commutes with~$D_N$. Then~$G$ also acts on~$M$ and~$E^+|_M$ such that the induced action on sections commutes with~$D_M$. One can define an {\em equivariant index} and an {\em equivariant~$\eta$-invariant} for all~$g \in G$ by
\begin{equation}\label{B1.2}\begin{aligned}
\ind_{\mathrm{APS},g}(D_N) &= \tr(g|_{\ker D_N^+}) - \tr(g|_{\ker D_N^-}) \\
\eta_{D_M,g}(s) &= \sum_{\lambda \in \Spec(D_M)\backslash\{0\}}\sign(\lambda) \cdot |\lambda|^{-s} \cdot \tr(g|_{E_\lambda}) \\
&= \frac{1}{\Gamma\bigl(\frac{s+1}{2}\bigr)}\int_0^\infty t^{\frac{s-1}{2}}\tr\Bigl(g\,D_M\,e^{-tD_M^2}\Bigr) \, dt\;.
\end{aligned}
\end{equation}
The equivariant~$\eta$-function has a meromorphic continuation to~$\C$, and 0 is a regular value. Again, we put~$\eta_g(D_M) = \eta_{D_M,g}(0)$. If one generalises the proof of Theorem~\ref{A1.T1} to this new setting, then the index density~$\alpha_0$ localises to the fix-point set~$N_g$ of~$g$. For a Dirac operator, we will write the equivariant index density as
$$ \hat{A}_g\bigl(TN,\nabla^{TN}\bigr)\, \ch_g\bigl(E/S,\nabla^E\bigr) \in \Omega^\bullet(N_g;o(N_g))\;,$$
where~$o(N_g)$ denotes the orientation line bundle. Note that~$\hat{A}_g(TN,\nabla^{TN})$ itself is a product of~$\hat{A}(TN_g,\nabla^{TN_g})$ and a contribution from the action of~$g$ on the normal bundle of~$N_g$ in~$N$.
Both forms are unique only up to a sign that depends on the choice of a lift of~$g$ to the spin group of~$N_g$, but their product is well-defined,
see the discussion in~\cite{BGV}.
We also put~$h_g(D_M)=\tr(g|_{\ker D_M})$.

\begin{Theorem}[Donnelly, \cite{Doeq}]\label{B1.T2}
The~$G$-equivariant index is given by
$$ \ind_{\mathrm{APS},g}(D_N) = \int_{N_g}\hat{A}_g\bigl(TN,\nabla^{TN}\bigr)\,\ch_g\bigl(E/S,\nabla^E\bigr) - \frac{\eta_g+h_g}2(D_M)\;.$$
\end{Theorem}

\begin{Remark}\label{B1.R1}
  The integral vanishes if~$g$ acts freely on~$N$, and the equivariant index is always a virtual character of~$G$. This has two consequences.
  \begin{enumerate}
  \item\label{B1.R1.1} There is an analogue of Corollary~\ref{A1.C3} with values in functions on~$G$ modulo virtual characters. For each~$g\in G$, the local contribution is an integral over~$M_g$. Hence, equivariant $\eta$-invariants are rigid modulo virtual characters for~$g\in G$ that act freely on~$M$.
  \item\label{B1.R1.2} Let~$D_M$ be a $G$-equivariant operator and let~$G_0\subset G$ denote a subset of elements that act freely on~$M$. If there is no virtual character~$\chi$ of~$G$ that extends~$\frac{\eta+h}2|_{G_0}$, and there is a compact manifold~$N$ with~$\partial N=M$ and~$D_N$ as in Donnelly's theorem, then some elements of~$g\in G$ will have fixpoints on~$N$.
  \end{enumerate}
\end{Remark}

If~$G$ is a compact connected Lie group, then the equivariant index theorem can be stated in a different way. Let~$\mathfrak g$ denote the Lie algebra of~$G$. Then consider Cartan's complex of {\em equivariant differential forms,\/}
$$ (\Omega_G^\bullet(N),d_{\mathfrak g}) = \Bigl((\Omega^\bullet(N)[\![\mathfrak g^*]\!])^G, d - \frac{\iota_X}{2\pi i}\Bigr)\;.$$
Here, a monomial in~$\mathfrak g^*$ of degree~$\ell$ with values in the~$k$-forms has total degree~$k + 2\ell$, and~$\iota_X$ denotes the inner product of a differential form with a variable Killing field~$X$, which is of total degree~$1 = -1+2$. Classical Chern-Weil theory generalises to~$G$-equivariant vector bundles with invariant connections, giving classes~$\hat{A}_X,\ch_X$ with values in the equivariant cohomology
$$H_G^\bullet(N;\R) = H^\bullet(\Omega_G^\bullet(N),d_{\mathfrak g})\;.$$
The classical equivariant index theorem can be stated in terms of these equivariant characteristic classes as explained by Berline, Getzler and Vergne~\cite{BGV}. 

Following Bismut's proof of the equivariant index theorem in~\cite{Binf}, put
$$ D_{X,t} = \sqrt{t}\,D_M + \frac{1}{4\sqrt{t}}\,c_X\;,$$
where~$c_X$ denotes Clifford multiplication with the Killing field associated to~$X \in \mathfrak g$. The {\em infinitesimally equivariant~$\eta$-invariant} of~$D_M$ is defined as
\begin{equation}\label{B1.3}
\eta_X(D_M) = \frac{2}{\sqrt{\pi}}\int_0^\infty \tr\Bigl(\frac{\partial D_{X,t}}{\partial t} e^{-D_{X,t}^2 - \mathcal L_X}\Bigr) \, dt \in \C[\![\mathfrak g^*]\!]\;,
\end{equation}
where~$\mathcal L_X$ denotes the Lie derivative. We can now state another version of the equivariant index theorem for manifolds with boundary.

\begin{Theorem}[\cite{Geta}]\label{B1.T3}
The equivariant index for~$g = e^{-X}$ is given by the formal power series
\begin{multline*}
  \ind_{\mathrm{APS},{e^{-X}}}(D_N) = \int_N \hat{A}_X\bigl(TN,\nabla^{TN}\bigr)\,\ch_X\bigl(E/S,\nabla^E\bigr)\\
  - \frac{\eta_X+h_{e^{-X}}}2(D_M) \in \C[\![\mathfrak g^*]\!]\;.
\end{multline*}
\end{Theorem}

This theorem can be deduced from the Bismut-Cheeger Theorem~\ref{B1.T1} by regarding fibre bundles with structure group~$G$ and applying the general Chern-Weil principle.

Another possible proof uses Donnelly's Theorem~\ref{B1.T2} and Bott's localisation formula.
Let~$\vartheta = \frac{1}{2\pi i} g^{TN}(\punkt,X)$ denote the dual of a variable Killing field, then
$$ d_X\Bigl(\frac{\vartheta_X}{d_X\vartheta_X}\alpha_X\Bigr) = \alpha_X - \frac{\alpha_X}{e_X(\nu)} \cdot \delta\;,$$
where~$\delta$ denotes the distribution of integration over the zero-set~$N_X$ of~$X$, and~$e_X(\nu)$ denotes the equivariant Euler class of the normal bundle~$\nu\to N_X$. Note that though the single terms are not defined on all of~$N$, the equation above still makes sense in an~$L^1$-sense, i.e., after integration over~$N$. Theorem~\ref{B1.T3} follows from Theorem~\ref{B1.T2} and the following result by an application of the localisation formula.  Both proofs are explained in~\cite{BGV} in the case~$\partial N=\emptyset$.

\begin{Theorem}[\cite{Geta}]\label{B1.T4} Assume that the Killing field~$X$ has no zeros on~$M$. Then
$$ \eta_X(D_M) = \eta_{e^{-X}}(D_M) + 2 \int_M \frac{\vartheta_X}{d_X\vartheta_X}\hat{A}_X\bigl(TM,\nabla^{TM}\bigr)\,\ch_X\bigl(E/S,\nabla^E\bigr) \in \C[\![\mathfrak g^*]\!]\;.$$
\end{Theorem}

One expects that for sufficiently small~$X \in \mathfrak g$, the formal power series in Theorems~\ref{B1.T3} and~\ref{B1.T4} converge, and that a similar formula also holds if~$X$ vanishes somewhere on~$M$. One also expects that one can apply both theorems to~$ge^{-X}$, where~$g \in G$ and~$X \in \mathfrak g$ with~$Ad_gX = X$, and the local contributions are integrated over~$N_g$. 

\begin{Remark}\label{B1.R2}
  In general, the equivariant $\eta$-invariant~$\eta_g(D_M)$ is only continuous in~$g$ as long as the fixpoint set~$M_g$ varies continuously in~$g$. In particular, it is usually singular at~$g=e$. The singularity near~$g=e$ is encoded in the integral in Theorem~\ref{B1.T4}. Arguing as in Remark~\ref{B1.R1}~\eqref{B1.R1.2}, we see that the singularity of the integral above at~$X=0$ contains information about fixpoints of elements of~$G$ on compact $G$-manifolds~$N$ with~$\partial N=M$.
\end{Remark}

\subsection{Homogeneous spaces}\label{B3}
It seems that the introduction of families and group actions is an unnecessary complication if one is mainly interested in the~$\eta$-invariants of section~\ref{A1}. In the presence of a Lie group action, Theorem~\ref{B1.T4} allows to split the infinitesimal $\eta$-invariant~$\eta_X(D)$ into a rigid global object~$\eta_{e^{-X}}(D)$ and a locally computable correction term if~$X\ne 0$ everywhere on~$M$. If both terms can computed, then their sum extends continuously to~$X=0$ and gives the ordinary $\eta$-invariant.

Assume that~$H\subset G$ are compact Lie groups, and let~$D$ be the geometric Dirac operator on the homogeneous space~$M=G/H$ with a normal metric. In~\cite{Gequi}, we consider the {\em reductive Dirac operator\/}~$\tilde D$. It is a selfadjoint differential operator with the same principal symbol as~$D$, but~$\tilde D$ is better adapted to homogeneous spaces. This operator was independently discovered by Kostant~\cite{Kostant}. If~$G$ and~$H$ are not of the same rank, then most elements~$g\in G$ act freely on~$G/H$. By Remark~\ref{B1.R1}~\eqref{B1.R1.1},
\begin{equation}\label{B3.1}
  \frac{\eta_g+h_g}2\bigl(\tilde D\bigr)-\frac{\eta_g+h_g}2(D)=\chi(g)
\end{equation}
for all~$g\in G$ that act freely on~$M$, where the {\em equivariant spectral flow\/}~$\chi$ is a virtual character of~$G$. Moreover~$\chi=0$ for the untwisted Dirac operator. On the other hand, the kernel of the reductive odd signature operator~$\tilde B$ has no topological significance, and hence the spectral flow does not vanish in general for~$D=B$.

Given three compact Lie groups~$H\subset K\subset G$, one considers the fibration~$G/H\to G/K$ with fibre~$K/H$. The equivariant $\eta$-invariant~$\eta_G(\tilde D)$ for~$G/H$ can be computed from the equivariant $\eta$-invariant of a reductive Dirac operator either on the base~$G/K$ or on the fibre~$K/H$, whichever is odd-dimensional. The formula is similar to the adiabatic limit formula in Theorem~\ref{B1.T1}, but an equivariant $\eta$-invariant appears instead of an $\eta$-form, and no limit has to be taken. Suppose that~$S\subset T$ are maximal tori of~$H$ and~$G$, we consider the fibrations~$G/S\to G/H$ and~$G/S\to G/T$. This way, the computation of~$\eta_G(\tilde D)$ is reduced in two steps to the computation of an equivariant $\eta$-invariant of a twisted Dirac operator on the flat torus~$T/S$, which vanishes unless~$\rk G-\rk H=\dim T-\dim S=1$. Although the formula for~$\eta_G(\tilde D)$ in~\cite{Gequi} contains representation theoretic expressions, explicit knowledge of the representations of~$G$ is not needed. In particular, the spectrum of~$\tilde D$ on~$M$ is not computed, in contrast to the examples in section~\ref{A2}.

In~\cite{Ghomogen}, a formula for the correction term in Theorem~\ref{B1.T4} is given, again using the fibrations~$G/S\to G/H$ and~$G/S\to G/T$ considered above. Combining this with~\eqref{B3.1}, one obtains a formula for~$\frac{\eta_X+h_{e^{-X}}}2(D)$ up to a virtual character of~$G$. Evaluating at~$X=0$ gives~$\frac{\eta+h}2(D)\in\R/\Z$. By estimation of sufficiently many small eigenvalues of~$D$ and~$\tilde D$, one can even determine the equivariant spectral flow. This method is applied to compute the Eells-Kuiper invariant of the Berger space~$SO(5)/SO(3)$ in~\cite{GKS}, see section~\ref{C5}.

\subsection{Orbifolds with boundary}\label{B4}
Theorem~\ref{A1.T1} is generalised to orbifolds by Farsi~\cite{Farsi}. For simplicity, we state the version for Dirac operators. Farsi's original theorem holds
in the generality of Theorem~\ref{A1.T1}. 

Let~$M$ be an~$n$-dimensional orbifold. In particular, for each~$p\in M$ there exists a local parametrisation of the form~$\psi \colon V \to \Gamma_p\backslash V \cong U \subset M$. Here, the {\em isotropy group\/}~$\Gamma_p\subset O(n)$ of~$p$ is a finite subgroup acting linearly on~$V \subset \R^n$.
If~$\gamma \in \Gamma$, let~$(\gamma)$ denote its conjugacy class, and let~$C_\Gamma(\gamma)$ denote its centraliser in~$\Gamma$. Then the {\em inertia orbifold}~$\Lambda M$ consists of all pairs~$(p,(\gamma))$ with~$(\gamma)$ a conjugacy class in~$\Gamma_p$. A parametrisation of~$\Lambda M$ around~$(p,(\gamma))$ is given by
$$ \psi_{(\gamma)} \colon C_\Gamma(\gamma)\backslash V^\gamma \to \psi(V^\gamma) \times \{(\gamma)\} \subset \Lambda M\;.$$
In general, the inertia orbifold is not effective. The multiplicity~$m(p,(\gamma))$ defines a locally constant function on~$M$ that says how many elements of the isotropy group~$C_\Gamma(\gamma)$ act trivially on the fixpoint set~$V^\gamma$.

An orbifold vector bundle~$E$ over~$M$ is given by trivialisations of~$\psi^*E \to V$ for all parametrisations~$\psi$, together with an action of the isotropy group~$\Gamma$ on~$\psi^*E$ and compatible gluing data. A smooth section is represented locally by a~$\Gamma$-equivariant section of~$\psi^*E$. There are natural notions of Dirac bundles and Dirac operators. Because~$\Gamma(E)$ is a vector space, one can define the index and the~$\eta$-invariant of a Dirac operator as before.

On~$\Lambda M$, one defines characteristic differential forms~$\hat{A}_{\Lambda M}$ and~$\ch_{\Lambda M}$ such that
\begin{align*}
\psi^*_{(\gamma)}\hat{A}_{\Lambda M}\bigl(TM,\nabla^{TM}\bigr) &= \frac{1}{m(\gamma)}\hat{A}_\gamma\bigl(TM,\nabla^{TM}\bigr) \in \Omega^\bullet(V^\gamma,o(V^\gamma)) \\
\text{and}\qquad \psi^*_{(\gamma)}\ch_{\Lambda M}\bigl(E/S,\nabla^E\bigr) &= \ch_\gamma\bigl(E/S,\nabla^E\bigr) \in \Omega^\bullet(V^\gamma)\;.
\end{align*}
Apart from the multiplicity, these forms are the same as in
Theorem~\ref{B1.T2}. In particular, the signs of both forms depend on the choice of a lift of~$\gamma$, but their product is well-defined.
The integrand~$\hat{A}_{\Lambda M}\bigl(TM,\nabla^{TM}\bigr)\wedge\ch_{\Lambda M}\bigl(E/S,\nabla^E\bigr)$ on~$\Lambda M$ is the same as in Kawasaki's index theorem, and on the regular part of~$M \cong M \times \{\id\} \subset \Lambda M$, it agrees with the classical index density on a manifold.

We now assume that~$N$ is an orbifold with boundary~$M$, and that~$D_N$, $D_M$ are Dirac operators satisfying~\eqref{A1.1}.

\begin{Theorem}[Farsi, \cite{Farsi}]\label{B1.T5}
The orbifold index under~$APS$ boundary conditions is given by
$$ \ind_{\mathrm{APS}}(D_N) = \int_{\Lambda N} \hat{A}_{\Lambda N}\bigl(TN,\nabla^{TN}\bigr)\, \ch_{\Lambda N}\bigl(E/S,\nabla^E\bigr) - \frac{\eta+h}{2}(D_M)\;.$$
\end{Theorem}

If~$N$ is a quotient of a compact manifold by a finite group of isometries, then Theorem~\ref{B1.T5} can be deduced from Theorem~\ref{B1.T2}. In general, one combines the proof of Kawasaki's index theorem with the proof of Theorem~\ref{A1.T1}.

\section{Properties of \texorpdfstring{$\eta$}{eta}-invariants}\label{P}
We state some formulas that do not directly follow from the Atiyah-Patodi-Singer index theorem and its generalisations in the previous section. The formulas are useful to understand properties of secondary invariants derived from $\eta$-invariants as in section~\ref{C}, and sometimes even to compute them.

\subsection{The adiabatic limit}\label{P1}
Let~$p \colon M \to B$ be proper Riemannian submersion with fibre~$F$ and~$TM=T^HM\oplus TF$ with~$T^HM = TF^\perp \cong p^*TB$. Write~$g^{TM} = g^{TF} \oplus p^*g^{TB}$ and define
$$ g_\eps^{TM} = g^{TF} \oplus \frac{1}{\eps^2} g^{TB}\;.$$
The limit~$\eps \to 0$ is called the {\em adiabatic limit}. As the distance between different fibres becomes arbitrarily large in the adiabatic limit, heat kernels of adapted Laplacians localise to a fibrewise operators as~$\eps \to 0$ for bounded times. This allows to localise a large part of the integral~\eqref{A1.3} to the fibres of~$p$.

Let~$(D_{M,\eps})_{\eps > 0}$ be a family of Dirac operators on a bundle~$E \to M$ that are compatible with the metrics~$g_\eps^{TM}$. We assume that the connections~$\nabla^{E,\eps}$ converge to a limit connection~$\nabla^{E,0}$. Associated to the limit~$\eps \to 0$, there exists a family of superconnections~$(\A_t)_{t > 0}$ as in section~\ref{B1}. The vertical Dirac operator~$D_F$ appears as the degree zero component of~$\A_1$. We assume that~$H = \ker D_F$ forms a vector bundle over~$B$. Then we can define the $\eta$-form~$\eta(\A) \in \Omega^\bullet(B)$ as in~\eqref{B1.1}. More precisely, if~$SB\to B$ is a local spinor bundle on~$B$, then there exists a fibrewise Dirac bundle~$W\to M$ such that~$E=p^*SB\otimes W$, and we consider the $\eta$-form of a superconnection~$\mathbb A$ on~$p_*W$.

The bundle~$H \to B$ with the connection induced by~$\nabla^{E,0}$ becomes a Dirac bundle on~$(B,g^{TB})$, and one can construct a limit Dirac operator~$D_B^0$ acting on~$H$. We assume that~$D_{M,\eps}$ can be continued analytically in~$\eps$ to~$\eps = 0$. Then~$\ker D_{M,\eps}$ has constant dimension for all~$\eps \in (0,\eps_0)$ if~$\eps_0 > 0$ is sufficiently small. There are finitely many {\em very small eigenvalues}~$\lambda = \lambda_\nu(\eps)$ of~$D_{M,\eps}$ such that
$$ \lambda_\nu(\eps) = O(\eps^2)\qquad \text{and}\qquad 0 \not= \lambda_\nu(\eps)\, \text{for}\, \eps \in (0,\eps_0)\;.$$

\begin{Theorem}[Bismut-Cheeger,\cite{BCeta}; Dai, \cite{Dai}]\label{P1.T1}
Under the assumptions above and for~$\eps \in (0,\eps_0)$, one has
$$ \lim_{\eps \to 0} \eta(D_{M,\eps}) = \int_B \hat{A}(TB,\nabla^{TB}) 2\eta(\A) + \eta(D_B^0) + \sum_\nu \sign(\lambda_\nu(\eps))\;.$$
\end{Theorem}

Both the Levi-Civita connection~$\nabla^{TM,\eps}$ and the connection~$\nabla^{E,\eps}$ converge as~$\eps \to 0$, so one can still define Chern-Simons classes as in~\eqref{A1.4}. Moreover, the spectral flow of~\eqref{A1.7} eventually becomes constant by our assumptions above, so one can recover~$\eta(D_M)$ from
\begin{multline*}
\frac{\eta+h}{2}(D_M) = \lim_{\eps \to 0}\Bigl(\frac{\eta+h}{2}(D_{M,\eps}) - \SF\bigl((D_{M,s})_{s \in (\eps,1]}\bigr)\Bigr) \\
+ \int_M \Bigl(\tilde{\hat A}\bigl(TM,\nabla^{TM,0},\nabla^{TM}\bigr)\,\ch\bigl(E/S,\nabla^{E,0}\bigr) \\
+ \hat{A}\bigl(TM,\nabla^{TM}\bigr)\,\tilde{\ch}\bigl(E/S,\nabla^{E,0},\nabla^E\bigr)\Bigr)\;.
\end{multline*}

Theorem~\ref{P1.T1} can be generalised to Seifert fibrations. Here, a {\em Seifert fibration} is a map~$p \colon M \to B$, where~$M$ is a manifold and~$B$ an orbifold, such that locally for a parametrisation~$\psi$ as in section~\ref{B4}, $p$ pulls back to
$$ \psi^*p \colon \psi^*M \cong V \times F \to V\;.$$
Then we call~$F$ the generic fibre of~$p$. Equivalently, a Seifert fibration is a Riemannian foliation of~$M$ with compact leaves. We define metrics~$g_\eps^{TM}$ as above. 

Over the inertia orbifold~$\Lambda B$, we define an equivariant~$\eta$-form~$\eta_{\Lambda B}(\A)$ such that
$$ \psi^*_{(\tilde\gamma)}\eta_{\Lambda B}(\A) = \frac{1}{\sqrt{\pi}}\int_0^\infty \tr\biggl(\tilde\gamma \frac{\partial \A_t}{\partial t} e^{-\A_t^2}\biggr) \, dt \in \Omega^\bullet(V^\gamma)\;.$$
Again, the sign of~$\eta_{\Lambda B}(\A)$ depends on the choice of a certain lift~$\tilde\gamma$ of~$\gamma$, but the integrand~$\hat{A}_{\Lambda B}\bigl(TB,\nabla^{TB}\bigr)\, 2\eta_{\Lambda B}(\A)$ in the theorem below is well-defined.

We assume that~$\ker D_F$ forms an orbifold vector bundle over~$B$ and define~$\eta(D_B^H)$ as in section~\ref{B4}.

\begin{Theorem}[\cite{Gorbi}]\label{P1.T2}
Under the assumptions above and for~$\eps \in (0,\eps_0)$, one has
$$ \lim_{\eps \to 0} \eta(D_{M,\eps}) = \int_{\Lambda B}\hat{A}_{\Lambda B}\bigl(TB,\nabla^{TB}\bigr)\, 2\eta_{\Lambda B}(\A) + \eta\bigl(D_B^H\bigr) + \sum_\nu\sign(\lambda_\nu(\eps))\;.$$
\end{Theorem}
It is likely that this result still holds if~$M$ is an orbifold, provided the generic fibres are still compact manifolds.

\begin{Remark}\label{P1.R1} In principle, Theorems~\ref{P1.T1} and~\ref{P1.T2} simplify the computations of~$\eta$-invariants and other invariants derived from them as in section~\ref{C}. However, the~$\eta$-forms needed are at least as difficult to compute as the~$\eta$-invariants of the fibres. There are explicit formulas for circle bundles in~\cite{Zhang} and three-sphere bundles in~\cite{Geta}. In~\cite{DZ}, the Kreck-Stolz invariants of section~\ref{C4} are computed this way for circle bundles. And in section~\ref{B3}, we have exhibited a method of computation if the structure group is compact and the fibre a quotient of compact Lie groups.

If the family~$M \to B$ bounds a family~$N \to B$ as in section~\ref{B1}, one can use the original Atiyah-Patodi-Singer Theorem~\ref{A1.T1} in place of Theorem~\ref{P1.T1}, see~\cite{BS} for the case of circle bundles. Similarly, one can use Theorem~\ref{B1.T5} in place of Theorem~\ref{P1.T2} if the generic fibre~$F$ bounds a compact manifold and one can construct an orbifold fibre bundle~$N \to B$ that bounds~$M$. Nevertheless, the computation of the local index density still requires some work.
\end{Remark}

\subsection{Gluing Formulas}\label{P2}

The~$\eta$-invariants of connected sums can be computed by applying the APS Index Theorem~\ref{A1.T1} to the boundary connected sum~$M_1 \times [0,1]\mathbin{\natural}M_2 \times [0,1]$ with boundary~$-(M_1 \# M_2) \sqcup M_1 \sqcup M_2$. Because Pontrijagin forms are conformally invariant, one can choose the geometry in such a way that the index density vanishes completely. Hence,
$$ \frac{\eta+h}{2}(D_1 \# D_2) = \frac{\eta+h}{2}(D_1) + \frac{\eta+h}{2}(D_2) \in \R/\Z$$
under suitable geometric assumptions. As a consequence, many of the invariants introduced in section~\ref{C} are additive under connected sums. We will now describe the behaviour of $\eta$-invariants under gluing along more complicated hypersurfaces.

We assume that~$M$ can be cut along a hypersurface~$N$ in two pieces~$M_1$ and~$M_2$. We also assume that~$N$ has a neighbourhood~$U$ isometric to~$N \times (-\eps,\eps)$. Let~$A = D_M$ be a Dirac operator on~$M$ that is of a form similar to~\eqref{A1.1} on~$U$, with~$B = D_N$ a Dirac operator on~$N$ and~$\nu = c_M\bigl(\frac{\partial}{\partial t}\bigr)$. If~$D_N$ is invertible, one can define~$\eta$-invariants for the operators~$D_{M_i} = D_M|_{M_i}$ under APS boundary conditions similar to~\eqref{A1.2}.

\begin{Theorem}[Wojciechowski, \cite{Woj1}; Bunke, \cite{Buglu}]\label{P2.T1}
If~$D_N$ is invertible, then
$$ \eta(D_M) = \eta(D_{M_1}) + \eta(D_{M_2})\qquad \in \R/\Z\;.$$
\end{Theorem}

This formula holds in~$\R$ up to an integer correction term that is also described in~\cite{Buglu} and~\cite{Woj1}. If~$D_N$ is not invertible, one chooses Lagrangian subspaces~$L_1$, $L_2 \subset \ker(D_N)$, with respect to a symplectic structure on~$\ker(N)$ defined in terms of the Clifford volume element on~$N$. The APS boundary conditions modified by the projections onto these subspaces give rise to selfadjoint operators~$D_{M_i,L_i}$. Their~$\eta$-invariants are described by Lesch and Wojciechowski in~\cite{LeWoj}. Bunke and Wojciechowski generalise Theorem~\ref{P2.T1} to this setting in~\cite{Buglu} and~\cite{Woj2}. Their formula involves the Maslov index of the Lagrangians~$L_1$, $L_2$.

Gluing results for~$\eta$-invariants as in Theorem~\ref{P2.T1} allow to understand the behaviour of the secondary invariants of section~\ref{C} under operations like surgery. But since manifolds with boundary appear only as intermediate steps in these constructions, it would be nice to have a general gluing formula where no manifolds with boundary occur. Bunke states such a result in~\cite{Buglu}. 

\subsection{Embeddings}\label{P3}

In this section, let~$\iota \colon M \to N$ be a smooth embedding of compact spin manifolds. If~$D_M$ is a Dirac operator on~$M$, one constructs a~$K$-theoretic direct image~$D_N$ on~$N$ and compares the associated~$\eta$-invariants. Hence, the main result of this section is similar in spirit to Theorem~\ref{P1.T1} of Bismut-Cheeger and Dai.

More precisely, let~$SM \to M$ and~$SN \to N$ denote spinor bundles on~$M$ and~$N$. Then the normal bundle~$\nu \to M$ of the embedding has a spinor bundle~$S\nu \to M$ such that~$SN|_M \cong SM \otimes S\nu$. A {\em direct image} of a complex vector bundle~$V \to M$ consists of a complex vector bundle~$W = W^+ \oplus W^- \to N$ and a selfadjoint endomorphism~$A = a + a^*$ of~$W$ with~$a \colon W^+ \to W^-$, such that~$A$ is invertible on~$N \backslash M$ and degenerates linearly along~$M$, and one has an isomorphism~$\ker A \cong S\nu \otimes V$ that relates the compression of~$dA|_\nu$ to~$\ker A$ to Clifford multiplication by normal vectors. We also assume that~$V$, $W$ and~$\nu$ carry compatible metrics and connections, see~\cite{BZeta} for details.

Let~$\delta_M$ denote the current of integration on~$M \subset N$. Then there exists a natural current~$\gamma(W,\nabla^W,A)$ on~$N$ such that
$$ d\gamma(W,\nabla^W,A) = \ch(W^+,\nabla^{W^+}) - \ch(W^-,\nabla^{W^-}) - \hat{A}^{-1}(\nu,\nabla^\nu)\ch(V,\nabla^V) \delta_M\;.$$

\begin{Theorem}[Bismut-Zhang, \cite{BZeta}]\label{P3.T1} Under the assumptions above,
\begin{align*}
& \frac{\eta+h}{2}(D_N^{W^+}) - \frac{\eta+h}{2}(D_N^{W^-}) = \frac{\eta+h}{2}(D_M^V) + \int_N \hat{A}(TN,\nabla^{TN})\gamma(W,\nabla^W,A) \\
& \qquad + \int_M \tilde{\hat{A}}\bigl(TN|_M, \nabla^{TM \oplus \nu},\nabla^{TN}\bigr) \hat{A}^{-1}(\nu,\nabla^\nu)\ch(V,\nabla^V) \in \R/\Z\;.
\end{align*}
\end{Theorem}

One can get rid of the last term on the right hand side by assuming that~$M$ is totally geodesic in~$N$. Moreover, if~$A$ is~$\nabla^W$-parallel outside a small neighbourhood of~$M$ in~$N$, then~$\gamma(W,\nabla^W,A)$ is supported near~$M$. In this case, the difference of the~$\eta$-invariants of~$D_N^{W^+}$ and~$D_N^{W^-}$ localises near~$M$. On the other hand, let~$M = \emptyset$, so~$a\colon W^+ \to W^-$ is an isomorphism of vector bundles. Then Theorem~\ref{P3.T1} reduces to Corollary~\ref{A1.C3}~\eqref{A1.C3.1} with~$g_s$ constant.

\begin{Remark}\label{P3.R1}
  Theorem~\ref{P3.T1} is formally similar to Theorem~\ref{P1.T1}. In fact, since every proper map~$F\colon M\to N$ can be decomposed into the embedding~$M\to M\times N$ given by the graph of~$F$, followed by projection onto~$N$, both theorem can be combined to compute $\eta$-invariants of direct images under proper maps. These direct images carry additional geometric information (like a connection), so one would like to have a generalisation of topological $K$-theory that takes care of the relevant additional data. In principle, some kind of smooth $K$-theory should be the right choice for this, but it seems difficult to construct a smooth $K$-theory that covers both proper submersions and embeddings.
\end{Remark}

\section{Differential topological invariants}\label{C}
We have seen that~$\eta$-invariants have local variation formulas with respect to variations of the geometric structure, but they still contain global differential-topological information. The common theme of the following sections will be the construction of invariants that do not depend on the geometric structure of the manifold.

\subsection{Invariants of flat vector bundles}\label{C1}
Let~$\alpha \colon \pi_1(M) \to U(n)$ be a unitary representation of the fundamental group, then~$F_\alpha = \tilde{M} \times_\alpha \C^n \to M$ is a flat Hermitian vector bundle with holonomy~$\alpha$. In particular, we may regard the twisted odd signature operator~$B_\alpha$ acting on even~$F_\alpha$-valued smooth forms on~$M$. Because~$\ker(B_\alpha) = H^\ev(M;F_\alpha)$ is independent of the Riemannian metric on~$M$, the variation formula in Corollary~\ref{A1.C3} becomes
$$ \eta(B_\alpha^1)-\eta(B_\alpha^0) = \int_M n\tilde{L}\bigl(TM,\nabla^{TM,0},\nabla^{TM,1}\bigr)\in\R\;.$$
Note that~$\alpha$ enters on the right hand side only through its rank~$n = \ch(F_\alpha)$.

\begin{Theorem}[Atiyah-Patodi-Singer, \cite{APS2}] \label{C1.T1}
The~$\rho$-{\em invariant}
$$ \rho_\alpha(M) = \eta(B_\alpha) - n\, \eta(B) \in \R$$
is a diffeomorphism invariant of~$M$ and~$\alpha$.
\end{Theorem}

If~$\alpha$ factors through a finite group~$G$, one can consider the compact manifold~$\bar{M} = \tilde{M}/\ker \alpha$. Then~$G$ acts on~$\bar{M}$ with quotient~$M$, and one can compute~$\rho_\alpha(M)$ from the equivariant signature~$\eta$-invariants~$\eta_g(\bar{M})$ of~$M$. This proves in particular that~$\rho_\alpha(M)$ is rational in this case. The equivariant~$\eta$-invariants here are related to the invariants~$\sigma_g(M)$ considered in~\cite{AS3}.

If~$\pi_1(M)$ is torsion free and a certain Baum-Connes assembly map is an isomorphism, then~$\rho_\alpha(M)$ is a homotopy invariant.
This is proved by Keswani in~\cite{Keswani} and generalised by Piazza and Schick~\cite{PS2},
earlier similar results are due to Neumann~\cite{Neumann}, Mathai~\cite{Mathai} and Weinberger~\cite{Weinberger}.
Hence for such fundamental groups, $\rho_\alpha(M)$ behaves almost as a primary invariant.

If one replaces the odd signature operator in the construction of~$\rho_\alpha(M)$ by a different Dirac operator~$D$ on~$M$, one gets similar invariants with values in~$\R/\Z$ due to the possible spectral flow. However, instead of diffeomorphism invariants, one now obtains cobordism invariants.

\begin{Theorem}[Atiyah-Patodi-Singer, \cite{APS2}, \cite{APS3}]\label{C1.T2}
The~$\xi$-{\em invariant}
$$ \xi_\alpha(D) = \frac{\eta + h}{2}(D_\alpha) - n \cdot \frac{\eta + h}{2}(D) \in \R/\Z$$
is a cobordism invariant in the sense that~$\xi_\alpha(D) = 0$ if there exists a compact manifold~$N$ with~$M = \partial N$ such that~$D$ extends to an operator on~$N$ in the sense of~\eqref{A1.1} and~$\alpha$ extends to a representation of~$\pi_1(N)$.

The representation~$\alpha$ defines a class~$[\alpha] \in K^{-1}(M;\R/\Z)$ and the symbol of~$D$ gives~$\sigma \in K^1(TM)$. Then there exists a topological index~$\Ind_{[\alpha]}(\sigma) \in \R/\Z$, and
$$ \xi_\alpha(D_M) = \Ind_{[\alpha]}(\sigma)\;.$$
\end{Theorem}

\begin{Example}\label{C1.X0}
If~$M$ is spin and~$SM$ is a fixed spinor bundle on~$M$, then all other spinor bundles arise by twisting~$SM$ with real line bundles. In particular, the difference of~$\frac{\eta+h}2(D)$ for different spin structures is a $\xi$-invariant. Real line bundle are classified by~$H^1(M;\Z/2\Z)$. Dahl investigates these $\xi$-invariants for spin structures induced by the mod~$2$ reduction of integer classes, and also their dependence on the initial spin structure in~\cite{Dahl}.
\end{Example}

It is possible to define~$\xi_\alpha(D) \in \C/\Z$ for flat vector bundles associated to non-unitary representations~$\alpha \colon \pi_1(M) \to Gl(n,\C)$, see~\cite{APS3}, \cite{JW}. In this case, the imaginary part of~$\xi_\alpha(D)$ is related to the Kamber-Tondeur classes (also known as Borel classes) of~$\alpha$. Choose a Hermitian metric on~$F_\alpha$ and a unitary connection~$\nabla^u$ on~$F_\alpha$ and let~$D_u$ be the Dirac operator twisted by~$(F_\alpha,\nabla^u)$. Arguing as in~\cite{BL},
\begin{align*}
  \Im\xi_\alpha(D)
  &= \Im\biggl(\frac{\eta+h}{2}(D_\alpha) - \frac{\eta+h}{2}(D_u)\biggr)\\
  &= \Bigl(\hat{A}(TM)\, \Im \tilde{\ch}\bigl(F_\alpha,\nabla^u,\nabla^\alpha\bigr)\Bigr)[M]\;,
\end{align*}
and~$\Im \tilde{\ch}\bigl(F_\alpha,\nabla^u,\nabla^\alpha\bigr)\in H^\odd(M;\R)$ represents the Kamber-Tondeur class.

Assume that~$M$ is an~$m$-dimensional homology sphere, then the fundamental group~$\Gamma = \pi_1(M)$ satisfies~$\Gamma = [\Gamma,\Gamma]$. Let~$F_\alpha \to M$ be associated to a representation~$\alpha \colon \Gamma \to Gl(n,\C)$, classified by a map~$M \to BGL(n,\C)$. Then Quillen's plus construction by functoriality gives an element~$[M,\alpha]$ of the algebraic~$K$-group~$K_m(\C) = \pi_m(BGL(n,\C)^+)$ by
$$ S^m = M^+ \longrightarrow BGL(n,\C)^+\;.$$
If~$m$ is odd, clearly~$\hat{A}(TM) = 1 \in H^\ev(M;\Q)$ because~$M$ is a homology sphere, in particular, $\Im \xi_\alpha(D)=\Im \tilde{\ch}\bigl(F_\alpha,\nabla^u,\nabla^\alpha\bigr)[M]$ then gives the Borel regulator of~$[M,\alpha] \in K_m(\C)$, see~\cite{JW}. Jones and Westbury prove that the map~$(M,\alpha) \mapsto \xi_\alpha(D)$ induces an isomorphism~$K_1(\C) \cong \C/\Z$ for~$m = 1$, and an isomorphism of the torsion subgroup of~$K_m(\C)$ with~$\Q/\Z$ for~$m > 1$ odd.

\begin{Example}\label{C1.X1} Let~$M = \Gamma\backslash SL(2,\C)/SU(2)$ be a hyperbolic homology three-sphere, and let~$\alpha \colon \Gamma \to SL(2,\C)$ be the representation corresponding to the embedding of~$\Gamma$ as a cocompact subgroup. By~\cite{JW},
$$ \Im \xi_\alpha(D) = -\frac{1}{4\pi^2}\operatorname{vol}(M)\;,$$
which proves that~$[M,\alpha]$ is never torsion.

Jones and Westbury also show that all torsion elements of~$K_3(\C)$ can be realised as~$[M,\alpha]$ where~$M$ now is a Seifert fibred three-manifold. A similar analysis of~$K_3(\R)$ is done in~\cite{Cm2}.
\end{Example}

The classical Lichnerowicz theorem asserts that a spin Dirac operator on a closed spin manifold~$N$ has vanishing index if~$N$ carries a metric of positive scalar curvature~$\kappa > 0$. It is shown in~\cite{APS2} that Lichnerowicz' theorem extends to compact spin manifolds~$N$ with totally geodesic boundary~$M = \partial N$. If~$N$ has~$\kappa > 0$, then so does~$M$, so~$\ind(D_N) = h(D_M) = 0$ for a spin Dirac operator~$D_N$ on~$N$, and Theorem~\ref{A1.T1} becomes
$$ \int_N\hat{A}(TN,\nabla^{TN}) = \frac{1}{2} \eta(D_M)\;.$$
The analogous statement also holds for the Dirac operator~$D_\alpha$ twisted by a flat vector bundle on~$N$ associated to~$\alpha \colon \pi_1(N) \to U(n)$.

Let us call two closed Riemannian spin manifolds~$(M_0,g_0)$, $(M_1,g_1)$ {\em spin$^+$-cobordant} if there exists a compact Riemannian spin manifold~$(N,g)$ with totally geodesic boundary~$(M_1,g_1) - (M_0,g_0)$. If~$M_0 = M = M_1$ and there exists a family~$(g_t)_{t \in [0,1]}$ of positive scalar curvature metrics, then the metric~$g_{\varphi(t/a)}\oplus dt^2$ on~$N = M \times [0,a]$ has~$\kappa > 0$ for a sufficiently large, where~$\varphi \colon [0,1] \to [0,1]$ is smooth and locally constant 0 (1) near 0 (1).  Thus metrics in the same connected component of the moduli space of positive scalar curvature metrics are~spin$^+$-cobordant.

\begin{Theorem}[Atiyah-Patodi-Singer, \cite{APS2}; Botvinnik-Gilkey \cite{BoGi}]\label{C1.T3}
The number
$$ \bar{\xi}_\alpha(M,g) = \eta(D_\alpha) - n \cdot \eta(D) \in \R$$
is a~spin$^+$-cobordism invariant in the sense that~$\bar{\xi}_{\alpha_0}(M_0,g_0) = \bar{\xi}_{\alpha_1}(M_1,g_1)$ if~$(N,g)$ is a spin$^+$-cobordism of~$(M_0,g_0)$ and~$(M_1,g_1)$ and~$\alpha_0,\alpha_1$ extend to a representation of~$\pi_1(N)$.
\end{Theorem}

This result is used by Botvinnik and Gilkey to construct and detect Riemannian metrics with~$\kappa > 0$ lying in countably many different connected components in the moduli space of positive scalar curvature metrics on~$M$, whenever~$M$ has a non trivial finite fundamental group and admits at least one metric of positive scalar curvature~\cite{BoGi}, \cite{BoGi2}. Apart from the computation of~$\bar{\xi}_\alpha$ for sufficiently many examples, the proof relies on the surgery techniques for positive scalar curvature metrics introduced by Gromov and Lawson in~\cite{GL}. Using various generalisations of~$\eta$-invariants, the results of Botvinnik-Gilkey are extended to manifolds~$M$ whose fundamental group contains torsion by Leichtnam and Piazza~\cite{LePi} and Piazza and Schick~\cite{PiSch}. Other (and in fact earlier) results in this direction will be discussed at the end of section~\ref{C3}.

On the other hand, if~$\pi_1(M)$ is torsion free and a certain Baum-Connes assembly map is an isomorphism, then~$\bar\xi_\alpha(D)=0$ for the untwisted Dirac operator~$D$ by a result of Piazza and Schick~\cite{PS2}. Hence, for such fundamental groups, the invariant~$\bar\xi_\alpha(D)$ behaves similar as the index of the untwisted Dirac operator in Lichnerowicz' theorem.

\subsection{The Adams \texorpdfstring{$e$}{e}-invariant}\label{C2}
In this subsection, we regard a framed bordism invariant. Recall that a closed manifold~$M$ is framed by an embedding~$M \hookrightarrow \R^n$ for~$n$ sufficiently large together with a trivialisation of the normal bundle~$\nu \to M$ of the embedding. This defines a stable parallelism of~$TM$, i.e., a trivialisation of~$TM \oplus \R^{n-m}$, because
$$ TM \oplus \R^{n-m} \cong TM \oplus \nu \cong M \times \R^n\;.$$
In fact, framings and stable parallelisms are equivalent notions. By the Pontrijagin-Thom construction, the framed bordism classes of manifolds of dimension~$m$ are in bijection with the~$m$-th stable homotopy group~$\pi_m^s$ of spheres. 

Let~$M$ be a framed closed manifold of dimension~$4k - 1$ with parallelism~$\pi$. Then~$M$ carries a preferred spin structure. Because the spin cobordism group in dimension~$4k-1$ is trivial, there exists a compact spin manifold~$N$ with~$\partial N = M$. Because~$TM$ is stably trivial, there exists a well-defined relative class~$\hat{A}(TN) \in H^\bullet(N,M;\Q)$ and one defines
\begin{equation}\label{C2.1}
e(M,\pi) = \begin{cases} \hat{A}(TN)[N] \qquad & \text{if $k$ is even, and} \\
\frac{1}{2} \hat{A}(TN)[N] & \text{if $k$ is odd}\;. \end{cases}
\end{equation}

On the other hand, since~$TM$ is stably trivial, we can consider the Chern-Simons class~$\tilde{\hat A}(TM,\nabla^{TM},\nabla^\pi)$, where~$\nabla^{TM}$ denotes the connection induced on~$TM \oplus \R^{n-m}$ by the Levi-Civita connection with respect to a Riemannian metric, and~$\nabla^\pi$ the connection induced by the trivialisation. It follows from Corollary~\ref{A1.C3} that
$$ \frac{\eta + h}{2} (D_M) + \int_M \tilde{\hat A}\bigl(TM,\nabla^{TM},\nabla^\pi\bigr)$$
is invariant under variation of~$g$ modulo~$\Z$ if~$k$ is even, and modulo~$2\Z$ if~$k$ is odd due to a quaternionic structure on the spinor bundle of~$M$.

\begin{Theorem}[Atiyah-Patodi-Singer, \cite{APS2}]\label{C2.T1}
The $e$-invariant of a framed~$4k-1$-dimensional manifold~$(M,\pi)$ is given by
$$ e(M,\pi) = \eps(k)\biggl(\frac{\eta + h}{2}(D_M) + \int_M \tilde{\hat A}\bigl(TM,\nabla^{TM},\nabla^\pi\bigr)\biggr)\in\Q/\Z$$
with~$\eps(k) = 1$ if~$k$ is even and~$\eps(k) = \frac{1}{2}$ if~$k$ is odd.
\end{Theorem}

Seade uses this formula to determine the $e$-invariants of quotients of~$S^3$ in~\cite{Seade}.

\begin{Example}\label{C2.X2}
Let~$H(n) \subset Gl_{n+2}(\R)$ denote the~$2n+1$-dimensional Heisenberg group and let~$\Gamma(n) = H(n) \cap Gl_{n+2}(\Z)$ be the subgroup with integer entries. Then~$TH(n)$ is trivialised by right translation, and this descends to a trivialisation~$\pi$ of~$TH(n)/\Gamma(n)$. For odd~$n = 2k-1$, the~$e$-invariant is calculated by Deninger and Singhof in~\cite{DS},
$$ e(H(n)/\Gamma(n),\pi) = -(-1)^k \eps(k)\zeta(-n) + \delta(n)\;,$$
where~$\delta(1) = \frac{1}{2}$ and~$\delta(n) = 0$ otherwise. Here~$\zeta$ is the Riemann zeta function. Comparing with the possible values of~$e(M,\pi)$, one sees that~$e(H(n)/\Gamma(n))$ is a generator of~$\im(e \colon \pi_{4k-1}^s \to \Q/\Z)$ for odd~$k$ and twice a generator for even~$k$.

For the proof, the Dirac operator~$D$ is replaced by~$\tilde{D}$, where~$\tilde{D} - D$ is an operator of order 0. The spectrum and the~$\eta$-invariant of~$\tilde{D}$ are computed explicitly.  Since~$e(H(n)/\Gamma(n),\pi) - \eps(k) \cdot \frac{\eta+h}{2}(\tilde{D}) \in \R/\Z$ is given as the integral of a locally defined invariant density on~$H(n)/\Gamma(n)$, an argument involving finite covering spaces allows to reconstruct the~$e$-invariant from the~$\eta$-invariant of the modified operator~$\tilde{D}$.
\end{Example}

Bunke and Naumann give a similar description of the $f$-in\-var\-i\-ant~\cite{BN}. Their formula uses $\eta$-invariants on manifolds with boundary that are related to a certain elliptic genus.

\subsection{The Eells-Kuiper invariant}\label{C3}
In this section, we consider closed oriented spin manifolds~$M$ of dimension~$m = 4k-1$ such that
\begin{equation}\label{C3.1}
H^{4l}(M;\R) = 0 \qquad \text{for all} \, l \ge 1\;.
\end{equation}
If~$M$ bounds a compact spin manifold~$N$,
this conditions allows to define relative Pontrijagin classes~$p_j(TN)\in H^{4j}(N,M;\Q)$ for~$1\le j<k$.
We express the universal characteristic classes~$\hat{A}$ and~$L$ in terms of Pontrijagin classes as in~\eqref{A1.5},~\eqref{A1.6}. Then there exists a unique constant~$t_k \in \Q$ such that the homogeneous component~$(\hat{A} - t_kL)^{[4k]}$ in degree~$4k$ does not involve~$p_k$. With~$\eps(k)$ as in Theorem~\ref{C2.T1}, the Eells-Kuiper invariant of~$M$ is defined in~\cite{EK} as
\begin{equation}\label{C3.0}
  \mu(M)=\eps(k)\,\bigl(t_k\,\sign(N)+(\hat{A} - t_kL)(TN)[N,M]\bigr)\in\Q/Z\;.
\end{equation}

Condition~\eqref{C3.1} allows one to express the Pontrijagin forms~$p_j(TM,\nabla^{TM})$ with respect to some Riemannian metric~$g$ on~$M$ as
\begin{equation}\label{C3.2}
  p_j(TM,\nabla^{TM}) = d\hat{p}_j(TM,\nabla^{TM})
\end{equation}
for~$1 \le j < k$. Moreover, $\hat{p}_j(TM,\nabla^{TM}) \in \Omega^{4j-1}(M)/\im d$ is unique because~$H^{4j-1}(M;\R) = 0$ by Poincar\'e duality.
Replacing one factor~$p_j$ in each monomial of~$(\hat{A} - t_kL)^{[4k]}(TM)$ by~$\hat{p}_j$, we obtain a natural class~$\alpha(TM,\nabla^{TM}) \in H^{4k-1}(M;\R) = \Omega^{4k-1}(M)/\im d$ such that
$$ \alpha\bigl(TM,\nabla^{TM,1}\bigr) - \alpha\bigl(TM,\nabla^{TM,0}\bigr) = \Bigl(\tilde{\hat A} - t_k \tilde{L}\Bigr)\bigl(TM,\nabla^{TM,0}, \nabla^{TM,1}\bigr)$$
for any two connections~$\nabla^{TM,0},\nabla^{TM,1}$ on~$TM$.
Note that~$\alpha$ does not depend on the choice of the factors~$p_j$ above, because
$$(p_i \, \hat{p}_j - \hat{p}_i\,p_j)\bigl(TM,\nabla^{TM}\bigr) =  d\Bigl((\hat{p}_i \,\hat{p}_j)\bigl(TM,\nabla^{TM}\bigr)\Bigr)\;.$$
Let~$D$ again be the spin Dirac operator and~$B$ the odd signature operator on~$M$. The following result is a consequence of Theorems~\ref{A1.T1} and~\ref{A1.T2}.

\begin{Theorem}[Donnelly, \cite{Doek}; Kreck-Stolz, \cite{KS1}]\label{C3.T0}
The Eells-Kuiper invariant of~$M$ equals
\begin{equation*}
\mu(M) = \eps(k) \cdot \Bigl(\frac{\eta+h}{2}(D) - t_k\,\eta(B) - \alpha\bigl(TM,\nabla^{TM}\bigr)[M]\Bigr) \in \Q/\Z\;.
\end{equation*}
\end{Theorem}

\begin{Remark}\label{C3.R1}
Other interesting invariants have expressions similar to~\eqref{C3.0}.
\begin{enumerate}
\item\label{C3.R1.1}The Eells-Kuiper invariant distinguishes all diffeomorphism types of exotic spheres that bound parallelisable manifolds in dimension~$4k-1$ for~$k=1$, $2$, $3$, see~\cite{EK}. Stolz constructs a similar invariant that detects all exotic spheres bounding parallelisable manifolds in all dimensions~$4k-1$ in~\cite{Stolz}. Stolz' invariant also has a presentation in terms of $\eta$-invariants and Cheeger-Simons correction terms.
\item\label{C3.R1.2} Rokhlin's theorem says that the signature of a spin manifold in dimension~$8k+4$ is divisible by~$16$. One can define a secondary Rokhlin number in~$\R/16\Z$ for spin structures on $8k+3$-dimensional manifolds. Lee and Miller express the Rokhlin number as a linear combination of $\eta$-invariants as above and without local correction terms~\cite{ML}. In particular, condition~\eqref{C3.1} is not needed. If spin structures differ only by the mod~$2$ reduction of an integer cohomology class, then Dahl proved that the Rokhlin number mod~$8$ remains unchanged~\cite{Dahl}.
\end{enumerate}
\end{Remark}

Note that~$\mu(M)$ changes sign if the orientation of~$M$ is reversed, and that~$\mu$ is additive under connected sums. Hence, given any closed spin~$4k-1$-manifold~$M$ that satisfies assumption~\eqref{C3.1} and an exotic sphere~$\Sigma$, $\mu(M \# \Sigma^{\# r})$ takes as many different values in~$\Q/\Z$ as~$\mu(\Sigma^{\# r})$ does for~$r \in \Z$. This way, one can construct and detect a certain number of exotic smooth structures on manifolds~$M$ for which~$\mu(M)$ is defined.

As an example, for~$k = 2$ the Eells-Kuiper invariant becomes
$$ \mu(M) = \frac{\eta+h}{2}(D) + \frac{\eta}{2^5 \cdot 7}(B) - \frac{1}{2^7 \cdot 7}(p_1\,\hat{p}_1)\bigl(TM,\nabla^{TM}\bigr)[M]\;.$$
This invariant is one of the main ingredients in the diffeomorphism classification of~$S^3$-bundles over~$S^4$ by Crowley-Escher~\cite{CE}, and also in the examples discussed in section~\ref{C5}.

We now come back to manifolds of positive scalar curvature. Note that~$\mu(M)$ is not a spin-cobordism invariant because~$\mu(M_1) - \mu(M_0)$ depends on~$t_k \cdot \sign(N)$ if~$\partial N = M_1 - M_0$ by~\eqref{C3.0}. We thus cannot expect to refine~$\mu(M)$ to a spin$^+$-cobordism invariant. Thus, we call two positive scalar curvature metrics~$g_0$, $g_1$ on~$M$ {\em concordant} if there exists a positive scalar curvature metric~$g$ on~$M \times [0,T]$ for some~$T > 0$ that is isometric to~$g_0 \times dt^2$ on~$M \times [0,\eps)$ and to~$g_1 \times dt^2$ on~$M \times (T-\eps,T]$. From the discussion preceding Theorem~\ref{C1.T3}, we see that metrics in the same connected component of the space of scalar curvature metrics are concordant. 

\begin{Theorem}[Kreck-Stolz, \cite{KS3}]\label{C3.T1}
Let~$M$ be a closed spin~$4k-1$-manifold satisfying~\eqref{C3.1}. Then the refined Eells-Kuiper invariant
$$ \bar{\mu}(M,[g]) = \eps(k) \Bigl(\frac{\eta}{2}(D) - t_k \, \eta(B) - \alpha\bigl(TM,\nabla^{TM}\bigr)[M]\Bigr) \in \R$$
is well-defined on concordance classes~$[g]$ of positive scalar curvature metrics on~$M$. Moreover, if~$[g_0], [g_1]$ are two such concordance classes, then
$$ \bar{\mu}(M,[g_1]) - \bar{\mu}(M,[g_0]) \in \Z\;.$$
\end{Theorem}

All members of the family of Aloff-Wallach spaces~$M\cong SU(3)/S^1$ allow metrics of positive sectional curvature, and the numbers~$\bar\mu(M)$ are computed in~\cite{KS2}.
For Witten's family of Ricci-positive homogeneous Einstein manifolds~$M \cong SU(3) \times SU(2) \times S^1/SU(2) \times S^1 \times S^1$, the invariants~$\bar\mu(M)$ are computed in~\cite{KS1}.

\begin{Theorem}[Kreck-Stolz, \cite{KS3}]\label{C3.T2}
\begin{enumerate}
\item\label{C3.T2.1} There exist closed manifolds with a non-connected moduli space of positive sectional curvature metrics.
\item\label{C3.T2.2} There exist closed manifolds for which the moduli space of Ricci positive metrics has infinitely many connected components.
\end{enumerate}
\end{Theorem}

\subsection{Kreck-Stolz invariants of complex and quaternionic line bundles}\label{C4}
In~\cite{KS1}, Kreck and Stolz define three invariants that determine the diffeomorphism type of certain 7-manifolds completely. For spin manifolds, their first invariant is precisely the Eells-Kuiper invariant of the previous section. The other two invariants use Dirac operators twisted by complex line bundles. For non-spin-manifold, similar invariants are defined that use a spin$^c$-Dirac operator in place of the Dirac operator. We will restrict attention to the spin case for simplicity. 

Thus assume that~$M$ is a closed simply connected spin 7-manifold with
\begin{equation}\label{C4.-2}
  H^1(M) = H^3(M) = 0, \qquad H^2(M) \cong \Z, \qquad \text{and}\qquad H^4(M) \cong \Z/\ell \Z\;,
\end{equation}
where~$H^4(M)$ is generated by the square of a generator of~$H^2(M)$. In particular, condition~\eqref{C3.1} holds. Since~$H^2(M)$ classifies complex line bundles, for each class~$a \in H^2(M)$, there exists a complex line bundle~$L \to M$ with Chern class~$c_1(L)=a$, which is unique up to isomorphism. Let~$\nabla^L$ be a connection, then as in~\eqref{C3.2} above, there exists a unique class~$v(L,\nabla^l) \in \Omega^3(M)/\im d$ such that
$$ dv(L,\nabla^L) = c_1(L,\nabla^L)^2\;.$$
We define a universal formal power series~$\ch'$ in~$c_1$ such that
$$ \ch(L) - 1 - c_1(L) = c_1^2(L)\, \ch'(L)\;.$$
Now, let~$D^L$ denote the Dirac operator twisted by~$(L,\nabla^L)$, and put
\begin{multline}\label{C4.-1}
s_M(a)= \frac{\eta+h}{2}(D^L) - \frac{\eta+h}{2}(D) \\
 - \Bigl(\hat{A}\bigl(TM,\nabla^{TM}\bigr)\,\bigl(v\,\ch'\bigr)\bigl(L,\nabla^L\bigr)\Bigr)[M] \in \Q/\Z\;.
\end{multline}
Let~$u\in H^2(M)$ be a generator, then the remaining two Kreck-Stolz invariants are given by
$$ s_2(M) = s_M(u) \qquad\text{and}\qquad s_3(M) = s_M(2u)\;.$$
Clearly, $s_2$ and~$s_3$ determine~$s_M$ completely. Also, one can recover the linking form on~$H^4(M)$ and the half Pontrijagin class~$\frac{p_1}{2}(TM) \in H^4(M)$ from~$s_2$ and~$s_3$. Indeed, $M$ bounds a compact spin manifold~$N$, and by Theorem~\ref{A1.T1},
\begin{equation}\label{C4.0}
s_M(a) = \bigl(\hat{A}(TM)(\ch(L)-1)\bigr)[N,M] = \biggl(\frac{a^2}{24}\Bigl(a^2 - \frac{p_1}{2}(TM)\Bigr)\biggr)[N,M]\;.
\end{equation}
In particular,
$$ 24s_M(a) = \lk_M\Bigl(a^2,a^2 - \frac{p_1}{2}(TM)\Bigr)\in\Q/\Z\;.$$

Hepworth generalises the Kreck-Stolz classification in his thesis~\cite{Hep} to simply connected closed spin 7-manifolds with
$$ H^1(M) = H^3(M) = 0,\qquad H^2(M) \cong \Z^r\qquad \text{and}\qquad \# H^4(M) < \infty\;,$$
such that~$H^4(M)$ is generated by products of elements of~$H^2(M)$ and~$\frac{p_1}{2}(TM)$.

A 7-manifold~$M$ is called {\em highly connected\/} if~$\pi_1(M) = \pi_2(M) = 0$. If~$\pi_3(M)$ is finite, then
$$ H^1(M) = H^2(M) = H^3(M) = 0\qquad\text{and}\qquad \# H^4(M) < \infty\;.$$
Since~$H^4(M)$ is not necessarily generated by~$\frac{p_1}{2}(TM)$, the results of Hepworth do not apply. Crowley has shown in~\cite{Crow} that a highly connected 7-manifolds is determined up to diffeomorphism by its Eells-Kuiper invariant and a quadratic form~$q_M \colon H^4(M) \times H^4(M) \to \Q/\Z$ satisfying
\begin{align*} q_M(a+b) &= q_M(a) + q_M(b) + \lk_M(a,b)\\
\text{and}\qquad q_M(-a) &= q_M(a) + \lk_M\Bigl(a,\frac{p_1}{2}(TM)\Bigr)\;.
\end{align*}
Note that these properties do not define~$q_M$ uniquely if~$H^4(M)$ has 2-torsion. An extrinsic definition of~$q_M$ using a handlebody~$N$ with~$\partial N = M$ can be found in~\cite{Crow}. The quadratic form~$q_M$ can also be recovered from a Kreck-Stolz type invariant~$t$ that we now describe.

Assume that~$M$ is a closed~$4k-1$-dimensional spin manifold satisfying
\begin{equation}\label{C4.1}
H^3(M;\R) = H^4(M;\R) = 0\;.
\end{equation}
Let~$H \to M$ be a quaternionic Hermitian line bundle. Equivalently, $H$ is a complex rank 2 vector with structure group~$SU(2)$. In particular, the determinant line bundle~$\det H$ is trivialised. Then the Chern character of~$H$ is a formal power series in~$c_2$, and there exists a formal power series~$\ch'$ in~$c_2$ such that
$$ 2 - \ch(H) = c_2(H) \cdot \ch'(H)$$ 
for all quaternionic line bundles~$H$. 

We fix a compatible connection~$\nabla^H$ on~$H \to M$. By assumption~\eqref{C4.1} and as in~\eqref{C3.2}, there exists a unique class~$\hat{c}_2(H,\nabla^H) \in \Omega^3(M)/\im d$ such that
$$ d\hat{c}_2(H,\nabla^H) = c_2(H,\nabla^H)\;.$$
Let~$D^H$ denote the Dirac operator twisted by~$H$ and note that~$S \otimes H$ carries a quaternionic structure if and only if~$S$ carries a real structure and vice versa. Let~$\eps(k)$ be as in Theorem~\ref{C2.T1}. In~\cite{CG}, we define the~$t$-{\em invariant} of~$H$ in analogy with~\eqref{C4.-1} by
\begin{multline}\label{C4.2}
t_M(H) = \eps(k+1) \biggl(\frac{\eta+h}{2}(D^H) - (\eta+h)(D)\\
+ \Bigl(\hat{A}\bigl(TM,\nabla^{TM}\bigr)\bigl(\hat{c}_2\,\ch'\bigr)\bigl(H,\nabla^H\bigr)\Bigr)[M]\biggr) \in \Q/\Z\;.
\end{multline}

If~$M$ is a highly connected closed 7-manifold with~$H^4(M)$ finite, then for each~$a \in H^4(M)$, there exists a quaternionic line bundle~$H \to M$ with~$c_2(H) = a$, and similar as in~\eqref{C4.0}, we find that
$$ q_M(a) = 12\,t_M(H)\;.$$

Note that the invariants~$t_M$ and~$s_M$ are related. Let~$L\to M$ be a complex line bundle with~$c_1(L)=a$. Then~$H=L\oplus\bar L$ carries a natural quaternionic structure. It follows from~\eqref{C4.-1} and~\eqref{C4.2} that~$t_M(H)=2\eps(k+1)\,s_M(a)$, so~$t_M$ generalises~$s_M$ in dimension~$8\ell-1$.

\begin{Example}\label{C4.X1}
Let~$\pi \colon M \to S^4$ be the unit sphere bundle of a real vector bundle~$W \to S^4$ of rank 4, and pick a quaternionic line bundle~$H \to S^4$, such that
$$ n = e(W),\qquad p = \frac{p_1}{2}(W),\qquad\text{and}\qquad a = c_2(H) \in \Z \cong H^4(S^4)\;.$$
Note that~$n$ and~$p$ are of the same parity. Then~$c_2(\pi^*H) = a \in \Z/n \Z \cong H^4(M)$ by the Gysin sequence. As shown in~\cite{CG},
$$ t_M(\pi^*H) = \frac{a(p-a)}{24n}\qquad\text{and}\qquad q_M(a) = \frac{a(p-a)}{2n} \in \Q/\Z\;.$$
Together with the computation of the Eells-Kuiper invariant
$$ \mu(M) = \frac{p^2-n}{2^5 \cdot 7 \cdot n} \in \Q/\Z$$
in~\cite{CE}, one can recover the Crowley-Escher diffeomorphism classification of~$S^3$-bundles over~$S^4$.
\end{Example}

The above example already shows that different quaternionic line bundles can have the same second Chern class, but different~$t$-invariants. In fact, the classifying space~$BSU(2) \cong\HH P^\infty$ for quaternionic line bundles is not a~$K(\pi,4)$ because
$$ \pi_{\ell + 1}(BSU(2)) \cong \pi_\ell(SU(2)) \cong \pi_\ell(S^3)$$
by the exact sequence of the fibre bundle~$ESU(2) \to BSU(2)$. Hence, $c_2$ alone does not classify quaternionic line bundles. Take~$M = S^7$ as an example, then quaternionic line bundles are classified by elements of
$$ \pi_7(BSU(2)) \cong \pi_6(S^3) \cong \Z/12\Z\;.$$
We prove in~\cite{CG} that on highly connected 7-manifold~$M$ with~$\pi_3(M)$ finite as above, $\pi_6(S^3)$ acts simply transitively on the set of isomorphism classes of quaternionic line bundles with a fixed second Chern class. The group~$\pi_6(S^3)$ acts freely by a clutching construction over a small~$S^6 \subset M$, and this action is detected by the~$t$-invariant. 

The~$t$-invariant also distinguishes all quaternionic line bundles on~$S^{11}$, but not on~$S^{15}$. Regard the sequence of Hopf fibrations and inclusions
$$\begin{array}{ccccccccccc}
\dots\! & & & &\! S^{4k-1}\! & & & &\! S^{4k+3}\! & & \\
& \searrow &  & \swarrow & & \searrow &  & \swarrow &  & \searrow &\\
& &\! \HH P^{k-1}\! & & & &\! \HH P^k\! &  &  &  &\! \dots
\end{array} \qquad \;.$$
Here, $\HH P^{k}$ decomposes along~$S^{4k-1}$ into a~$4k$-disk and a 4-disk bundle over~$\HH P^{k-1}$. A quaternionic line bundle on~$\HH P^{k-1}$ can be extended to~$\HH P^k$ if and only if its pullback to~$S^{4k-1}$ is trivial.
The~$t$-invariant on~$S^{4k-1}$ is thus an obstruction against such an extension. By cellular approximation, quaternionic line bundles on~$\HH P^k$ are classified by homotopy classes of maps~$\HH P^k \to \HH P^k \subset \HH P^\infty \cong BSU(2)$. If we compute the~$t$-invariants on~$S^7, S^{11},\dots,S^{4k-1}$ for a quaternionic line bundle with a given second Chern class, we recover precisely the obstructions against self maps of~$\HH P^k$ found by Feder and Gitler in~\cite{FeGi}.

Finally, Crowley also defines an analogous quadratic form~$q_M$ on~$H^8(M)$ for highly connected $15$-manifolds~$M$ in~\cite{Crow}. An intrinsic formula for~$q_M$ will probably involve the unique string structure on~$M$ in the same way that~\eqref{C4.2} above uses the unique spin structure.

\subsection{Seven-manifolds of positive curvature}\label{C5}
Riemannian metrics of positive sectional curvature on closed manifolds
are a rare phenomenon,
and sharp conditions for their existence are far from being understood.
Apart from the obvious symmetric examples, few other manifolds are known.
Many of these other examples are seven-dimensional manifolds that are either of Kreck-Stolz type~\eqref{C4.-2} or highly connected.
The homogeneous Aloff-Wallach spaces~$SU(3)/U(1)$ and their biquotient analogues, the Eschenburg space, have been classified using Kreck-Stolz invariants in~\cite{KS2}, \cite{AMP} and~\cite{Kruggel}. Kruggel uses a cobordism with lens spaces, whose $\eta$-invariants have already been given in~\cite{APS2}.

The Berger space~$SO(5)/SO(3)$ is diffeomorphic to a particular $S^3$-bundle over~$S^4$. For the proof in~\cite{GKS}, one needs to know that it is homeomorphic to such a bundle by~\cite{KiSh}. Then the Eells-Kuiper invariant of~$SO(5)/SO(3)$ together with the classification of all $S^3$-bundles over~$S^4$ in~\cite {CE} suffices to determine the diffeomorphism type.

One is still interested in finding new examples of positive curvature metrics.
Grove, Wilking and Ziller~\cite{GWZ} give two families~$(P_k)$, $(Q_k)$ of 7-manifolds
and one exceptional space~$R$,
which possibly allow such metrics
and contain new examples.
The spaces~$P_k$ are highly connected, whereas~$Q_k$ and~$R$ are of Kreck-Stolz type.
In~\cite{GVZ}, Grove, Verdiani and Ziller constructed a positive sectional
curvature metric on~$P_2$ (note that~$P_1=S^7$);
another construction is due to Dearricott~\cite{Dear}.
On the other hand, the space~$R$ does not carry a metric of cohomogeneity one with positive sectional curvature by a result of Verdiani and Ziller~\cite{VZ}.

The spaces~$P_k$ form Seifert fibrations with generic fibre~$S^3$ over
some base orbifold~$B_k$ as indicated in~\cite{GWZ}.
We apply Theorem~\ref{P1.T2} to determine the $\eta$-invariants in~\eqref{C3.1} and~\eqref{C4.2} in the adiabatic limit
and compute~$\mu(P_k)$ and~$t_{P_k}$ for all~$P_k$.

\begin{Theorem}[\cite{Gorbi}]\label{Thm2}
  The 
  Eells-Kuiper invariant of~$P_k$ is given by
  \begin{gather*}
    \mu(P_k)=-\frac{4k^3-7k+3}{2^5\cdot 3\cdot 7}\quad\in\quad\Q
    /\Z\;.
    \tag{1}\label{Thm2.1}
  \end{gather*}
  Crowley's quadratic form~$q$ on~$H^4(P_k)\cong\Z/k\Z$ is given by
  \begin{gather*}
    q(\ell)=\frac{\ell(\ell-k)}{2k}\quad\in\quad\Q/\Z\;.
    \tag{2}\label{Thm2.2}
  \end{gather*}
\end{Theorem}

By comparing these values with the corresponding values for $S^3$-bundles
over~$S^4$ in~\cite{CE} and~\cite{CG}, see Example~\ref{C4.X1},
one can construct manifolds that are diffeomorphic to~$P_k$.

\begin{Theorem}[\cite{Gorbi}]\label{Thm3}
  Let~$E_{k,k}\to S^4$ denote the principal $S^3$-bundle 
  with Euler class~$k\in H^4(S^4)\cong\Z$,
  and let~$\Sigma_7$ denote the exotic seven sphere
  with~$\mu(\Sigma_7)=\frac1{28}$.
  Then there exists an orientation preserving diffeomorphism
	$$P_k\cong E_{k,k}\mathbin{\#}\Sigma_7^{\#\frac{k-k^3}6}\;.$$
  In particular, $P_k$ and~$E_{k,k}$ are homeomorphic.
\end{Theorem}

This result also implies that~$P_2$ with reversed orientation is diffeomorphic to some $S^3$-bundle
over~$S^4$, and to~$T_1S^4\#\Sigma_7$, where~$T_1S^4$ denotes the unit tangent bundle
of~$S^4$.

\end{document}